\documentclass[12TP,draft]{article}

\begin{document}

\begin{center}
\LARGE\noindent\textbf{ On Hamiltonian Bypasses in one Class of Hamiltonian Digraphs }\\

\end{center}
\begin{center}
\noindent\textbf{Samvel Kh. Darbinyan and Iskandar A. Karapetyan}\\

Institute for Informatics and Automation Problems, Armenian National Academy of Sciences

E-mails: samdarbin@ipia.sci.am, isko@ipia.sci.am\\
\end{center}

\textbf{Abstract}

Let $D$ be a strongly connected directed graph of order $n\geq 4$ which satisfies the following condition (*): for every  pair of non-adjacent vertices $x, y$ with a common in-neighbour $d(x)+d(y)\geq 2n-1$ and $min \{ d(x), d(y)\}\geq n-1$. In \cite{[2]}  (J. of Graph Theory 22 (2) (1996) 181-187)) J. Bang-Jensen, G. Gutin and H. Li proved that $D$ is Hamiltonian. In [9] it was shown that if $D$  satisfies the condition (*) and the minimum semi-degree of $D$ at least two, then either $D$ contains a pre-Hamiltonian cycle (i.e., a cycle of length $n-1$) or $n$ is even and $D$ is isomorphic to the complete bipartite digraph (or to the complete bipartite digraph minus one arc) with partite sets of cardinalities of $n/2$ and $n/2$. In this paper we show that if the minimum out-degree of $D$ at least two and the minimum in-degree of $D$ at least three, then $D$ contains also a Hamiltonian bypass, (i.e., a subdigraph is obtained from a Hamiltonian cycle by reversing exactly one arc).\\ 


\section {Introduction} 

The directed graph (digraph) $D$ is Hamiltonian  if it contains a Hamiltonian cycle,
i.e., a cycle that includes every vertex of $D$. A Hamiltonian bypass in $D$ is a subdigraph obtained from a Hamiltonian cycle by reversing exactly one arc. We recall the following well-known degree conditions (Theorems 1-5) that guarantee that a digraph is Hamiltonian. 

 \textbf{Theorem 1} (Nash-Williams \cite{[14]}). Let $D$ be a digraph of order $n$ such that for every vertex $x$, $d^+(x)\geq n/2$ and $d^-(x)\geq n/2$, then $D$ is Hamiltonian.

 \textbf{Theorem 2} (Ghouila-Houri \cite{[12]}). Let $D$ be a strong digraph of order $n$. If $d(x)\geq n$ for all vertices $x\in V(D)$, then $D$ is Hamiltonian.

 \textbf{Theorem 3} (Woodall \cite{[16]}). Let $D$ be a digraph of order $n\geq 2$. If $d^+(x)+d^-(y)\geq n$ for all pairs of vertices $x$ and $y$ such that there is no arc from $x$ to $y$, then $D$ is Hamiltonian.
 
\textbf{Theorem 4} (Meyniel \cite{[13]}). Let $D$ be a strong digraph of order $n\geq 2$. If $d(x)+d(y)\geq 2n-1$ for all pairs of non-adjacent vertices in $D$, then $D$ is Hamiltonian.

It is easy to see that Meyniel's theorem is a common generalization of Ghouila-Houri's and Woodall's theorems. For a short proof of Theorem 1.3, see \cite{[5]}. \\

C. Thomassen \cite{[15]} (for $n=2k+1$) and S. Darbinyan \cite{[6]}
(for $n=2k$) proved the following:

\textbf{Theorem 5} \cite{[15], [6]}. If $D$ is a digraph of order
$n\geq 5$ with minimum degree at least $n-1$ and with minimum
semi-degree at least $n/2-1$, then $D$ is Hamiltonian (unless some
extremal cases which are characterized).\\

In view of the next theorems we need the following definitions.

\textbf{Definition 1}. Let $D_0$ denote any digraph of order $n\geq 5$, $n$ odd, such that $V(D_0)=A\cup B$, where $A\cap B=\emptyset$, $A$ is an independent set with $(n+1)/2$ vertices, $B$ is a set of $(n-1)/2$ vertices inducing any arbitrary subdigraph, and $D_0$ has $(n+1)(n-1)/2$ arcs between $A$ and $B$. Note that
$D_0$ has no Hamiltonian bypass. 

\textbf{Definition 2}. For any $k\in [1, n-2]$ let $D_1$ denote a digraph of order $n\geq 4$, obtained from $K^*_{n-k}$ and $K^*_{k+1}$ by identifying a vertex of the first with a vertex of the second. 
Note that $D_1$ has no Hamiltonian bypass.

\textbf{Definition 3}. By $T(5)$ we denote a tournament of order 5 with vertex set $V(T(5))=\{x_1, x_2, x_3, x_4, y\}$ and arc set $A(T(5))=\{x_ix_{i+1}/ i\in [1,3]\}\cup \{x_4x_1, x_1y, x_3y, yx_2, yx_4,x_1x_3, x_2x_4\}$. $T(5)$  has no Hamiltonian bypass.\\

In \cite{[4]} it was proved that if a digraph $D$ satisfies the condition of Nash-Williams' or Ghouila-Houri's or Woodall's theorem, then $D$ contains a Hamiltonian bypass. In [4] the following theorem was also proved:

 \textbf{Theorem 6} (Benhocine \cite{[4]}). Every strongly 2-connected digraph of order $n$ and minimum degree at least $n-1$ contains a Hamiltonian bypass, unless $D$ is isomorphic to a digraph of type $D_0$.\\

In \cite{[7]} the first author proved the following theorem:

 \textbf{Theorem 7} (Darbinyan \cite{[7]}). Let $D$ be a strong digraph of order $n\geq 3$. If $d(x)+d(y)\geq 2n-2$ for all pairs of non-adjacent vertices in $D$, then $D$ contains a Hamiltonian bypass unless it is isomorphic to a digraph of the set $D_0\cup \{D_1, T_5, C_3\}$, where $C_3$ is a directed cycle of length 3.\\

For $n\geq 3$ and $k\in [2,n]$, $D(n,k)$ denotes the digraph of order $n$ obtained from a directed cycle $C$ of length $n$ by reversing exactly $k-1$ consecutive arcs. The first author \cite{[7],[8]} has studied the problem of the existence of $D(n,3)$ in digraphs with the condition of Meyniel's theorem and in oriented graphs with large in-degrees and out-degrees.

 \textbf{Theorem 8} (Darbinyan \cite{[7]}). Let $D$ be a strong digraph of order $n\geq 4$. If $d(x)+d(y)\geq 2n-1$ for all pairs of non-adjacent vertices in $D$, then $D$ contains a $D(n,3)$.

\textbf{Theorem 9} (Darbinyan \cite{[8]}). Let $D$ be an oriented graph of order $n\geq 10$. If the minimum in-degree and out-degree of $D$ at least $(n-3)/2$, then $D$ contains a $D(n,3)$.\\

Each of  Theorems 1-5  imposes a degree condition on all pairs of non-adjacent vertices (or on all vertices). The following  theorem (as well as Theorems 13 and 14) imposes a degree condition only for some pairs of non-adjacent vertices. 

\textbf{Theorem 10} \cite{[2]} (Bang-Jensen, Gutin, H.Li \cite{[2]}). Let $D$ be a strong digraph of order $n\geq 2$. Suppose that 
$$min\{d(x),d(y)\}\geq n-1 \quad \hbox{and} \quad  d(x)+d(y)\geq 2n-1 \eqno (*) 
$$ 
for every pair of non-adjacent vertices $x,y$ with a common in-neighbour, then $D$ is Hamiltonian.\\

 In \cite{[9]} the following results were obtained: 

\textbf{Theorem 11} \cite{[9]}. Let $D$ be a strong digraph of order $n\geq 3$ with  the minimum semi-degree of $D$ at least two. Suppose that $D$ satisfies the condition (*). Then either $D$ contains a pre-Hamiltonian  cycle or $n$ is even and $D$ is isomorphic to the complete bipartite digraph or to the complete bipartite digraph minus one arc with partite sets of cardinalities $n/2$ and $n/2$. \\
 
In this paper using Theorem 11 we prove the following:\\

 \textbf{Theorem 12} (Main Result). Let $D$ be a strong digraph of order $n\geq 4$ with  the minimum out-degree at least two and with minimum in-degree at least three. Suppose that 
$$min\{d(x),d(y)\}\geq n-1 \quad \hbox{and} \quad  d(x)+d(y)\geq 2n-1 \eqno (*) 
$$ 
for every pair of non-adjacent vertices $x,y$ with a common in-neighbour. Then $D$ contains a Hamiltonian bypass.\\

\section {Terminology and Notations}

We shall assume that the reader is familiar with the standard
terminology on the directed graphs (digraph)
 and refer the reader to the  monograph of Bang-Jensen and Gutin \cite{[1]} for terminology not discussed here.
  In this paper we consider finite digraphs without loops and multiple arcs. For a digraph $D$, we denote
  by $V(D)$ the vertex set of $D$ and by  $A(D)$ the set of arcs in $D$. The order of $D$ is the number
  of its vertices. Often we will write $D$ instead of $A(D)$ and $V(D)$. The arc of a digraph $D$ directed from
   $x$ to $y$ is denoted by $xy$ or $x\rightarrow y$. If $x,y,z$ are distinct vertices in $D$, then $x\rightarrow y\rightarrow z$ denotes that $xy$ and $yz\in D$. Two distinct vertices $x$ and $y$ are adjacent if $xy\in A(D)$ or $yx\in A(D) $ (or both). By $a(x,y)$ we denote the number of arcs with end vertices $x$ and $y$, in particular, $a(x,y)$ means that the vertices $x$ and $y$ are non-adjacent. For disjoint subsets $A$ and  $B$ of $V(D)$  we define $A(A\rightarrow B)$ \,
   as the set $\{xy\in A(D) / x\in A, y\in B\}$ and $A(A,B)=A(A\rightarrow B)\cup A(B\rightarrow A)$. If $x\in V(D)$
   and $A=\{x\}$ we write $x$ instead of $\{x\}$. If $A$ and $B$ are two distinct subsets of $V(D)$ such that every vertex of  $A$ dominates every vertex of $B$, then we say that $A$ dominates $B$, denoted by $A\rightarrow B$. The out-neighborhood of a vertex $x$ is the set $N^+(x)=\{y\in V(D) / xy\in A(D)\}$ and $N^-(x)=\{y\in V(D) / yx\in A(D)\}$ is the in-neighborhood of $x$. Similarly, if $A\subseteq V(D)$, then $N^+(x,A)=\{y\in A / xy\in A(D)\}$ and $N^-(x,A)=\{y\in A / yx\in A(D)\}$. The out-degree of $x$ is $d^+(x)=|N^+(x)|$ and $d^-(x)=|N^-(x)|$ is the in-degree of $x$. Similarly, $d^+(x,A)=|N^+(x,A)|$ and $d^-(x,A)=|N^-(x,A)|$. The degree of the vertex $x$ in $D$ is defined as $d(x)=d^+(x)+d^-(x)$ (similarly, $d(x,A)=d^+(x,A)+d^-(x,A)$). The path (respectively, the cycle) consisting of the distinct vertices $x_1,x_2,\ldots ,x_m$ ( $m\geq 2 $) and the arcs $x_ix_{i+1}$, $i\in [1,m-1]$  (respectively, $x_ix_{i+1}$, $i\in [1,m-1]$, and $x_mx_1$), is denoted  by $x_1x_2\cdots x_m$ (respectively, $x_1x_2\cdots x_mx_1$). We say that $x_1x_2\cdots x_m$ is a path from $x_1$ to $x_m$ or is an $(x_1,x_m)$-path. For a cycle  $C_k:=x_1x_2\cdots x_kx_1$ of length $k$, the subscripts considered modulo $k$, i.e., $x_i=x_s$ for every $s$ and $i$ such that  $i\equiv s\, (\hbox {mod} \,k)$. If $P$ is a path containing a subpath from $x$ to $y$ we let $P[x,y]$ denote that subpath. Similarly, if $C$ is a cycle containing vertices $x$ and $y$, $C[x,y]$ denotes the subpath of $C$ from $x$ to $y$. A digraph $D$ is strongly connected (or, just, strong) if there exists a path from $x$ to $y$ and a path from $y$ to $x$ for every pair of distinct vertices $x,y$.
  For an undirected graph $G$, we denote by $G^*$ the symmetric digraph obtained from $G$ by replacing every edge $xy$ with the pair $xy$, $yx$ of arcs.  $K_{p,q}$ denotes the complete bipartite graph  with partite sets of cardinalities $p$ and $q$.   For integers $a$ and $b$, $a\leq b$, let $[a,b]$  denote the set of
all integers which are not less than $a$ and are not greater than
$b$.  By $D(n;2)=[x_1x_{n}; x_1x_2\ldots  x_n]$ is denoted the Hamiltonian bypass obtained from a Hamiltonian cycle $x_1x_2\ldots x_nx_1$ by reversing the arc $x_nx_1$.

\section { Preliminaries }

The following well-known simple Lemmas 1 and 2 are the basis of our results
and other theorems on directed cycles and paths in digraphs. They
will be used extensively in the proof of our result. \\

\textbf{Lemma 1} \cite{[11]}. Let $D$ be a digraph of order $n\geq 3$
 containing a
 cycle $C_m$, $m\in [2,n-1] $. Let $x$ be a vertex not contained in this cycle. If $d(x,C_m)\geq m+1$,
 then  $D$ contains a cycle $C_k$ for all  $k\in [2,m+1]$. \fbox \\\\

The following lemma is a slight modification of a lemma by Bondy and Thomassen \cite{[5]}.

\textbf{Lemma 2}. Let $D$ be a digraph of order  $n\geq 3$
 containing a
 path $P:=x_1x_2\ldots x_m$, $m\in [2,n-1]$ and let $x$ be a vertex not contained in this path.
  If one of the following conditions holds:

 (i) $d(x,P)\geq m+2$;

 (ii) $d(x,P)\geq m+1$ and $xx_1\notin D$ or $x_mx_1\notin D$;

 (iii) $d(x,P)\geq m$, $xx_1\notin D$ and $x_mx\notin D$,

then there is an  $i\in [1,m-1]$ such that
$x_ix,xx_{i+1}\in D$ (the arc $x_ix_{i+1}$ is a partner of $x$), i.e., $D$ contains a path $x_1x_2\ldots
x_ixx_{i+1}\ldots x_m$ of length $m$ (we say that  $x$ can be
inserted into $P$ or the path $x_1x_2\ldots x_ixx_{i+1}\ldots$ $ x_m$
is extended from $P$ with $x$). \fbox \\\\

\textbf{Definition 4} (\cite{[1]}, \cite{[2]}). Let $Q=y_1y_2\ldots y_s$ be a path in a digraph $D$ (possibly, $s=1$) and let $P=x_1x_2\ldots x_t$, $t\geq 2$, be a path in $D-V(Q)$. $Q$ has a partner on $P$ if there is an arc (the partner of $Q$) $x_ix_{i+1}$ such that $x_iy_1,y_sx_{i+1}\in D$. In this case the path $Q$ can be inserted into $P$ to give a new $(x_1,x_t)$-path with  vertex set $V(P)\cup V(Q)$. The path $Q$ has a collection of partners on $P$ if there are integers $i_1=1<i_2<\cdots <i_m=s+1$ such that, for every $k=2,3,\ldots , m$ the subpath $Q[y_{i_{k-1}},y_{i_k-1}]$ has a partner on $P$.\\

\textbf{Lemma 3} (\cite{[1]}, \cite{[2]}, Multi-Insertion Lemma).  Let $Q=y_1y_2\ldots y_s$ be a path in a digraph $D$ (possibly, $s=1$) and let $P=x_1x_2\ldots x_t$, $t\geq 2$, be a path in $D-V(Q)$. If $Q$ has a collection of partners on $P$, then there is an $(x_1,x_t)$-path with vertex set $V(P)\cup V(Q)$. \fbox \\\\

The following lemma is obvious.

\textbf{Lemma 4}. Let $D$ be a digraph of order $n\geq 3$ and let $C:=x_1x_2\ldots x_{n-1}x_1$ be an arbitrary cycle of length $n-1$ in $D$ and let $y$ be the vertex  not on $C$. If $D$ contains no Hamiltonian bypass, then

(i) $d^+(y,\{x_i, x_{i+1}\})\leq 1$ and   $d^-(y,\{x_i, x_{i+1}\})\leq 1$ for all $i\in [1,n-1]$;

(ii) $d^+(y)\leq (n-1)/2$, $d^-(y)\leq (n-1)/2$ and $d(y)\leq n-1$;

(iii) if $x_ky, yx_{k+1}\in D$, then $x_{i+1}x_i\notin D$ for all $x_i\not=x_k$. \fbox \\\\

Let $D$ be a digraph of order $n\geq 3$ and let $C_{n-1}$ be a cycle of length $n-1$  
in $D$. If for the vertex $y\notin C_{n-1}$, $d(y)\geq n$, then we say that $C_{n-1}$ is a good cycle. Notice that, by Lemma 4(ii), if a digraph $D$ contains a good cycle, then $D$ also contains a Hamiltonian bypass. \\

 We now  need to state and prove some general lemmas.\\

\textbf{Lemma 5}. Let $D$ be a digraph of order $n\geq 6$  with minimum semi-degree at least two satisfying the condition (*).   Let $C:=x_1x_2\ldots x_{n-1}x_1$ be an arbitrary cycle of length $n-1$ in $D$ and let $y$ be the vertex  not on $C$. Then 
for any $i\in [1,n-1]$ the following holds:

(i) If $yx_i\notin D$ and $x_{i-2}x_i\notin D$, then $x_i$ has a partner on $C[x_{i+1},x_{i-2}]$ or $d(x_{i})\geq n-1$.

(ii) If $yx_i\notin D$ and $d(x_i)\leq n-2$, then   $x_i$ has a partner on $C[x_{i+1},x_{i-2}]$ or there is a  vertex $x_k\in C[x_{i+1},x_{i-2}]$ such that $\{x_k,x_{k+1}, \ldots , x_{i-2}\}\rightarrow x_i$.
 
(iii) If $yx_i\in D$ , $x_{i-2}x_i\notin D$ and $d^-(x_i)\geq 3$, then $x_i$ has a partner on $C[x_{i+1},x_{i-1}]$ or $d(x_{i})\geq n-1$.

\textbf{Proof}. (i) The proof is by contradiction. Assume that $x_i$ has no partner on $C[x_{i+1},x_{i-2}]$ and $d(x_{i})\leq n-2$. Since $d^-(x_i,\{y,x_{i-2}\})=0$ and $d^-(x_i)\geq 2$, there is an $x_k\in C[x_{i+1},x_{i-3}]$ such that $x_kx_i\in D$. From $x_k\rightarrow \{x_{k+1},x_i\}$ , $d(x_i)\leq n-2$, $x_ix_{k+1}\notin D$ and the condition (*) it follows that $x_{k+1}x_i\in D$. By a similar argument we conclude that $x_{i-2}x_i\in D$, which is a contradiction.

For the proofs of (ii) and (iii) we can use precisely the same arguments as in the proof of (i). \fbox \\\\

\textbf{Lemma 6}. Let $D$ be a digraph of order $n\geq 6$  with minimum semi-degree at least two satisfying the condition (*). Let $C:=x_1x_2\ldots x_{n-1}x_1$ be an arbitrary cycle of length $n-1$ in $D$ and let $y$ be the vertex  not on $C$. Then 

(i) If for some $i\in [1,n-1]$,  $x_iy\in D$ and $x_{i+1}, y$ are non-adjacent or (ii) $a(x_i,y)=2$ or

(iii) $d(y)\geq n-1$, then $D$ contains a Hamiltonian bypass.\\

\textbf{Proof}. (i) Assume that (i) is not true. Without loss of generality, we assume that $x_{n-1}y\in D$, $d(y, \{x_1,x_2,\ldots , x_a\})=0$ and $x_{a+1}, y$ are non-adjacent, where $a\geq 1$. Then $x_1$ and $y$ is a dominated pair of non-adjacent vertices with a common in-neighbour $x_{n-1}$. Therefore, by condition (*), $d(y)\geq n-1$ and $d(x_1)\geq n-1$. On the other hand, using Lemma 4(i) we obtain that  $d(y)\leq n-a$ and hence, $a=1$ and $d(y)=n-1$. This together with the condition (*) implies that $d(x_1)\geq n$. If $yx_2\in D$, then $x_{n-1}yx_2x_3\ldots x_{n-2}x_{n-1}$ is a good cycle in $D$ and therefore, $D$ contains  a Hamiltonian bypass. Assume therefore that $x_2y\in  D$. Since $d(x_1)\geq n$, by Lemma 2, $x_1$ has a partner on the path $C[x_2,x_{n-1}]$, i.e, there is an $(x_2,y)$-Hamiltonian path which together with the arc $x_2y$ forms a Hamiltonian bypass, which is a contradiction and completes the proof of (i).

(ii) It follows immediately from Lemmas 6(i) and 4(i).

(iii) Suppose, on the contrary, that $d(y)\geq n-1$ and $D$ contains no Hamiltonian bypass as well as no good cycle. By Lemma 4(ii), $d(y)=n-1$. From Lemma 6(ii) it follows that $a(y,x_i)=1$ for all $i\in [1,n-1]$. Without loss of generality, we may assume that (by Lemma 4(i)) 
$$
N^+(y)=\{x_1,x_3,\ldots ,x_{n-2}\} \quad \hbox {and} \quad N^-(y)=\{x_2,x_4,\ldots ,x_{n-1}\}. \eqno (1)
$$
Notice that 

(2) for every vertex $x_i$, $x_ix_{i-1}\notin D$ and $x_i$ has no partner on the path $C[x_{i+1},x_{i-1}]$ (for otherwise, $D$ contains a Hamiltonian bypass).

Assume first that $x_1x_3\in D$. Then it is not difficult to show that $x_2, x_{n-1}$ are non-adjacent and $x_2x_4\notin D$. Indeed, by (1) if $x_{n-1}x_2\in D$, then $D(n,2)=[x_1x_2; x_1x_3x_4yx_5\ldots x_{n-1}x_2]$; if $x_2x_4\in D$, then 
 $D(n,2)=[x_2x_3; x_2x_4x_5\ldots x_{n-1}yx_1x_3]$; and if $x_2x_{n-1}\in D$, then $D(n,2)=[x_2x_{n-1}; x_2yx_1x_3x_4\ldots x_{n-1}]$, in each case we have a contradiction. Now, since $x_{n-1}x_2\notin D$, $yx_2\notin D$ and $x_2$ has no partner on $C[x_3,x_1]$,  Lemma 5(i) implies that $d(x_2)\geq n-1$. On the other hand, using  Lemma 2(ii), $a(x_2,x_{n-1})=0$, $x_2x_4\notin D$ and (2), we obtain 
$$
n-1\leq d(x_2)=d(x_2,\{y,x_1,x_3\})+d(x_2,C[x_4,x_{n-2}])\leq n-2,
$$ a contradiction.

Assume second that $x_1x_3\notin D$. By the symmetry of the vertices $x_{n-2}$ and $x_1$ (by (1)), we also may assume that $x_{n-2}x_1\notin D$. Since $x_1$ has no partner on $C[x_3,x_{n-2}]$, again using Lemma 2(iii) and (2) we obtain that 
$$
d(x_1)=d(x_1,\{x_{n-1},x_2,y\})+d(x_1,C[x_3,x_{n-2}])\leq n-2.
$$
Therefore, by condition (*), we have that $x_1$ is adjacent with $x_3$ and $x_{n-2}$, i.e., $x_3x_1, x_1x_{n-2}\in D$, since $y\rightarrow \{x_{n-2},x_1,x_3\}$. Now it is easy to see that $\{x_3,x_4,\ldots , x_{n-2}\}\rightarrow x_1$, which contradicts that $x_{n-2}x_1\notin D$. In each case we obtain a contradiction, and hence the proof of Lemma 6(iii) is completed. \fbox \\\\

 The following simple observation is of importance in the rest of the paper.

\textbf{Remark}. Let $D$ be a digraph of order $n\geq 6$  with minimum semi-degree at least two satisfying the condition (*). Let $C:=x_1x_2\ldots x_{n-1}x_1$ be an arbitrary cycle of length $n-1$ in $D$ and let $y$ be the vertex  not on $C$. If $D$ contains no Hamiltonian bypass, then

(i) There are two distinct vertices $x_k$ and $x_l$ such that $\{x_k,x_{k+1}\}\cap \{x_l,x_{l+1}\}=\emptyset$,  $x_k\rightarrow y\rightarrow x_{k+1}$ and $x_l\rightarrow y \rightarrow x_{l+1}$ (by Lemmas 4(i) and 5(i)).

(ii) $x_{i+1}x_i \notin D$ for all $i\in [1,n-1]$.

(iii) If $y\rightarrow \{x_{i-1},x_{i+1}\}$ or $\{x_{i-1},x_{i+1}\}\rightarrow y$, then $x_i$ has no partner on the path $C[x_{i+1},x_{i-1}]$.

(iv) If $x_{i+1}x_{i-1}\in D$, then $d(x_i)\leq n-2$ (by Remark (i) and Lemma 6(ii)).
 \fbox \\\\

 \textbf{Lemma 7}. Let $D$ be a digraph of order $n\geq 6$  with minimum semi-degree at least two satisfying the condition (*).  Let $C:=x_1x_2\ldots x_{n-1}x_1$ be an arbitrary cycle of length $n-1$ in $D$ and let  $y$ be the vertex not on $C$. Assume that $y\rightarrow \{x_2,x_{n-1}\}$, $x_1\rightarrow y$ and $d(y,\{x_3,x_{n-2}\})=0$. Then $D$ contains a Hamiltonian bypass. 

\textbf{Proof}. The proof is by contradiction. Assume that $D$ contains no Hamiltonian bypass. From Remark (i) and Lemmas 6(i), 4(i) it follows that for some $j\in [4,n-4]$, $x_j \rightarrow  y \rightarrow x_{j+1}$.

Now we show that $x_{n-2}x_1\notin D$. Assume that this is not the case. Then  $x_{n-2}x_1\in D$ and $d(x_{n-1})\leq n-2$ (by Remark (iv)). Then, since $y\rightarrow \{x_2,x_{n-1}\}$, the condition (*) implies that $x_2$ and $x_{n-1}$ are adjacent, i.e., $x_2x_{n-1}\in D$ or $x_{n-1}x_2\in D$. If $x_2x_{n-1}\in D$, then $D(n,2)=[x_2x_{n-1}; x_2x_3\ldots x_{n-2}x_1yx_{n-1}]$, and if    $x_{n-1}x_2\in D$, then $D(n,2)=[x_{n-1}x_{1}; x_{n-1}x_2x_3\ldots x_{j}yx_{j+1}\ldots x_{n-2}x_1]$. In both cases we have a Hamiltonian bypass, a contradiction. Therefore $x_{n-2}x_1\notin D$.

 Now, since $x_1$ has no partner on $C[x_3,x_{n-2}]$, by Lemma 5(i), $d(x_1)\geq n-1$. On the other hand, from $d(y,\{x_3,x_{n-2}\})=0$ , $d(y)\leq n-2$ and the condition (*) it follows that $x_1x_3\notin D$ and $x_1x_{n-2}\notin D$ (in particular, $a(x_1,x_{n-2})=0$). Now using Lemma 2(ii) and Remark (ii) we obtain 
$$
n-1\leq d(x_1)=d(x_1,\{y,x_2,x_{n-1}\})+d(x_1,C[x_3,x_{n-3}])\leq n-2, $$ 
which is a contradiction. Lemma 7 is proved. \fbox \\\\

\textbf{Lemma 8}. Let $D$ be a digraph of order $n\geq 6$ with minimum semi-degree at least two satisfying the condition (*).  Let $C:=x_1x_2\ldots x_{n-1}x_1$ be an arbitrary cycle of length $n-1$ in $D$ and let  $y$ be the vertex not on $C$. If  $d^-(y)\geq 3$ and $y$ is adjacent with four consecutive vertices of the cycle $C$, then $D$ contains a Hamiltonian bypass.

\textbf{Proof}. Suppose, on the contrary, that $D$ contains no Hamiltonian bypass and no good cycle. Using Lemmas 6(i) and 4(i), without loss of generality, we can assume that  $\{x_{n-1},x_2\}\rightarrow y$ and $y\rightarrow \{x_1,x_3\}$. By Remarks (ii) and (iii) we have 

(3) $x_ix_{i-1}\notin D$ for each $i\in [1,n-1]$ and $x_1$ (respectively, $x_2$) has no partner on the path $C[x_2,x_{n-1}]$ (respectively, $C[x_3,x_{1}]$).
 
 If $x_{n-2}x_1\notin D$ and $x_{1}x_3\notin D$, then  using Lemma 2(iii) and (3) we obtain  
$$
d(x_1)=d(x_1,\{y,x_2,x_{n-1}\})+d(x_1,C[x_3,x_{n-2}])\leq n-2.
$$ 
Therefore, by condition (*), the vertices  $x_1,x_3$ are adjacent, since $y\rightarrow \{x_1,x_3\}$ ($x_1$ and $x_3$ has a common in-neighbour $y$). This means that $x_3x_1\in D$. Since $x_1$ has no partner on $C[x_3,x_{n-2}]$, it follows from $x_3\rightarrow \{x_1,x_4\}$ and $d(x_1)\leq n-2$ that $x_4x_1\in D$. Similarly, we conclude that $x_{n-2}x_1\in D$ which contradicts the assumption that  $x_{n-2}x_1\notin D$. Assume therefore that
$$
x_{n-2}x_1\in D \quad \hbox{or} \quad x_{1}x_3\in D. \eqno (4)
$$ 

Now we prove that $d(x_1)\geq n-1$. Assume that this is not the case, that is $d(x_1)\leq n-2$. Then again by condition (*)   $x_1,x_3$ are adjacent because of $y\rightarrow \{x_1,x_3\}$. Therefore $x_3x_1\in D$ or  $x_1x_3\in D$. If $x_3x_1\in D$, then it is not difficult to show that  $\{x_3,x_4,\ldots , x_{n-2}\}\rightarrow x_1$, i.e., $d(x_1)\geq n-1$, a contradiction. Assume therefore that $x_3x_1\notin D$ and  $x_1x_3\in D$. Then $C':=x_1x_3x_4\ldots x_{n-1}yx_1$ is a cycle of length $n-1$ missing the vertex $x_2$. Then $d(x_2)\leq n-2$ (by Remark (iv)). Now, since $x_2$ has no partner on $C[x_3,x_1]$ (by (3)) and $d^-(x_2,\{x_i,x_{i+1}\})\leq 1$ for all $i\in [3,n-2]$ (by Lemma 4(i)), it follows that $d^-(x_2,\{y,x_3,x_4,\ldots , x_{n-2}\})=0$. Then $x_{n-1}x_2\in D$ because of $d^-(x_2)\geq 2$. From $d^-(y)\geq 3$ and Lemma 6(i) it follows that there is a vertex $x_j\in C[x_4,x_{n-3}]$ such that $x_j\rightarrow y\rightarrow x_{j+1}$. Therefore $D(n,2)= [x_1x_2;x_1x_3x_4\ldots x_jyx_{j+1}\ldots x_{n-1}x_2]$ is a Hamiltonian bypass, a contradiction. This contradiction proves that  $d(x_1)\geq n-1$. 

Notice that $x_{n-1}x_2\notin D$, by Remark (iv). From (4) it follows that the following two cases are possible: $x_1x_3\in D$ (Case 1) or $x_1x_3\notin D$ and $x_{n-2}x_1\in D$ (Case 2).

\textbf{Case 1.} $x_1x_3\in D$. Then $d(x_2)\leq n-2$ (by Remark (iv)). It is easy to see that $x_2x_4\notin D$ and $x_{n-1}x_2\notin D$ (if $x_{n-1}x_2\in D$, then $D$ has a cycle of length $n-1$ missing $x_1$, and hence $d(x_1)\leq n-2$ which contradicts that $d(x_1)\geq n-1$). Thus, we have a contradiction against Lemma 5(i), since $d(x_2)\leq n-2$, $x_{n-1}x_2\notin D$ and $x_2$ has no partner on $C[x_3,x_{n-1}]$ (by (3)).

 \textbf{Case 2.} $x_1x_3\notin D$ and $x_{n-2}x_1\in D$. It is easy to see that $x_{n-1}x_2\notin D$ and $x_{n-3}x_{n-1}\notin D$. If $yx_{n-2}\in D$, then $x_{n-1}$ has no partner on $C[x_1,x_{n-2}]$. This together with $d(x_{n-1})\leq n-2$,   and $x_{n-3}x_{n-1}\notin D$ contradicts Lemma 5(i). Assume therefore that $y$ and $x_{n-2}$ are non-adjacent. Then, since $d(y)\leq n-2$ and $x_2y\in D$, we have that $x_2x_{n-2}\notin D$.

Assume first that $x_2x_{n-1}\in D$. Then $x_{n-2}x_2\notin D$ (for otherwise, the arc $x_{n-2}x_{n-1}\in C[x_3,x_{n-1}]$ is a  partner of $x_2$ on $C[x_3,x_{n-1}]$, a contradiction against (3)). Therefore $x_2$ and $x_{n-2}$ are non-adjacent. Now we have  $x_ix_2\in D$, where $i\in [4,n-3]$ since $d^-(x_2)\geq 2 $ and $d^-(x_2, \{y,x_3,x_{n-2},x_{n-1}\})=0$. It is not difficult to see that $d(x_2)\geq n-1$ (Lemma 5(i)). Then by Remark (ii) and Lemma 2 we obtain   
$$
n-1\leq d(x_2)=d(x_2,\{y,x_1,x_3,x_{n-1}\})+d(x_2,C[x_4,x_{n-3}])\leq n-1,
$$ 
i.e., $d(x_2)=n-1$ and $d(x_2, C[x_4,x_{n-3}])=n-5$. By Lemma 2, $x_2x_4$ and $x_{n-3}x_2\in D$. From $d(x_2)=n-1$ and the condition (*) it follows that $d(x_{n-2})\geq n$, since $x_2$ and $x_{n-2}$ are non-adjacent and have a common in-neighbour $x_{n-3}$. If $x_1x_{n-2}\in D$, then $D(n,2)=[x_1x_2; x_1x_{n-2}x_{n-1}yx_3x_4\ldots x_{n-3}x_2]$, a contradiction. Assume therefore that $x_1x_{n-2}\notin D$. Now we consider the cycle $C':=x_{n-3}x_2x_{n-1}yx_3x_4\ldots x_{n-3}$ of length $n-2$ which does not contain the vertices $x_{n-2}$ and $x_1$. Since $d(x_{n-2})\geq n$ and $x_1x_{n-2}\notin D$ (i.e., $a(x_1,x_{n-2})=1$), then $d(x_{n-2}, C')\geq n-1$. Therefore, by Lemma 1, there is a cycle, say $C''$, of length $n-1$ missing the vertex $x_1$. Then, since $d(x_1,C'')\geq n-1$, by Lemma 4(ii) $D$ contains a Hamiltonian  bypass. 

Assume second that $x_2$ and $x_{n-1}$ are non-adjacent. Then, since $d(x_{n-1})\leq n-2$, the condition (*) implies that $x_{n-2}x_2\notin D$. Then by Remark (ii) and Lemma 2(ii), we have
$$
 d(x_2)=d(x_2,\{y,x_1,x_3\})+d(x_2,C[x_4,x_{n-2}])\leq n-2.
$$
 This contradicts Lemma 5(i) (because of (3)) and completes the proof of Lemma 8. \fbox \\\\

From Lemmas 6, 7 and 8 immediately  the following lemma follows:

\textbf{Lemma 9}. Let $D$ be a digraph of order $n\geq 6$  with minimum out-degree at least two and with minimum in-degree at least three satisfying the condition (*). Let $C:=x_1x_2\ldots x_{n-1}x_1$ be an arbitrary cycle of length $n-1$ in $D$ and let $y$ be the vertex not on $C$. If the vertex $y$ is adjacent with three consecutive vertices of the cycle $C$, then  $D$ contains a Hamiltonian bypass. \fbox \\\\

\textbf{Lemma 10}. Let $D$ be a digraph of order $n\geq 6$  with minimum out-degree at least two and with minimum in-degree at least three satisfying the condition (*). Let $C:=x_1x_2\ldots x_{n-1}x_1$ be an arbitrary cycle of length $n-1$ in $D$ and let $y$ be the vertex not on $C$. If $D$ contains no Hamiltonian bypass and $x_{i-1}x_{i+1}\in D$ for some $i\in [1,n-1]$, then $d(x_i,\{x_{i-2},x_{i+2}\})=0$.

\textbf{Proof}. The proof is by contradiction. Without loss of generality, we may assume that $D$ has no Hamiltonian bypass, $x_{n-1}x_2\in D$ and $a(x_1,x_{3})\geq 1$ or $a(x_1,x_{n-2})\geq 1$. If $a(x_1,x_{3})\geq 1$ (respectively, $a(x_1,x_{n-2})\geq 1$), then, since $y$ is not adjacent with three consecutive vertices of $C$, by Remark (i) there exists a vertex $x_k\in C[x_{3},x_{n-2}]$ (respectively, $x_k\in C[x_{2},x_{n-3}]$) such that $x_k\rightarrow y\rightarrow x_{k+1}$. It is not difficult to see that $C':=C[x_{2},x_{k}]yC[x_{k+1},x_{n-1}]x_{2}$ is a cycle of length $n-1$ missing the vertex $x_1$, and $x_1$ is adjacent with three consecutive vertices of $C'$, namely with $x_{n-1}, x_2,x_3$ (respectively, $x_{n-2}, x_{n-1},x_2$), which is a contradiction against Lemma 9. Lemma 10 is proved. \fbox \\\\

\textbf{Lemma 11}. Let $D$ be a digraph of order $n\geq 6$  with minimum out-degree at least two and with minimum in-degree at least three satisfying the condition (*). Let $C:=x_1x_2\ldots x_{n-1}x_1$ be an arbitrary cycle of length $n-1$ in $D$ and let $y$ be the vertex not on $C$. $D$ contains no Hamiltonian bypass, then $x_{i+1}x_{i-1}\notin D$ for all $i\in [1,n-1]$.

\textbf{Proof}. The proof is by contradiction. Without loss of generality, we may assume that $x_3x_1\in D$. 

Assume first that the vertex $x_2$ has a partner on $ C[x_4,x_{n-1}]$, i.e., there is an $x_j\in C[x_4,x_{n-2}]$ such that $x_j\rightarrow x_2\rightarrow x_{j+1}$. From $d^-(y)\geq 3$ and Lemma 6(i) it follows that there exists a vertex $x_{k}\in C[x_3,x_{n-2}]$ distinct from $x_j$ such that $x_k \rightarrow y\rightarrow x_{k+1}$. Therefore, if  $k\geq j+1$, then $D(n,2)=[x_{3}x_{1}; x_{3}x_4\ldots x_{j}x_2x_{j+1}\ldots x_kyx_{k+1}\ldots x_{n-1}x_1]$, and if $k\leq j-1$, then $D(n,2)=[x_{3}x_{1}; x_{3}x_4\ldots x_{k}yx_{k+1}$ $\ldots x_{j}x_2x_{j+1}\ldots x_{n-1}x_1]$, a contradiction. 

Assume second that $x_2$ has no partner on $C[x_4,x_{n-1}]$. Since $x_3x_1\in D$, Lemma 10 implies that $x_2x_4\notin D$ and $x_{n-1}x_2\notin D$. Now using Lemma 2(iii) and Remark (ii) we obtain  
$$
d(x_2)=d(x_2,\{y,x_1,x_{3}\})+d(x_2,C[x_4,x_{n-1}])\leq n-2.
$$ 
This together with the condition (*) implies that $d^-(x_2,C[x_3,x_{n-1}])=0$. Therefore $d^-(x_2)\leq 2$, which contradicts that 
$d^-(x_2)\geq 3$. Lemma 11 is proved. \fbox \\\\

\section {The proof of the main result}

 \textbf{Proof of Theorem 12.} By Theorem 11 the digraph $D$ contains a cycle of length $n-1$ or $n$ is even and $D$ is isomorphic to the complete bipartite digraph (or to the complete bipartite digraph minus one arc) with partite sets of cardinalities $n/2$ and $n/2$. If $n\leq 5$ or $D$ contains no cycle of length $n-1$, then it is not difficult to check that $D$ contains a Hamiltonian bypass. Assume therefore that $n\geq 6$, $D$ contains a cycle of length $n-1$ and has no Hamiltonian bypass. From Lemma 9 it follows that if $C$ is an arbitrary cycle of length $n-1$ in $D$ and the vertex $y$ is not on $C$, then there are not three consecutive vertices of $C$ which are adjacent with $y$. Let $C:=x_1x_2\ldots x_{n-1}x_1$ be an arbitrary cycle of length $n-1$ in $D$ and let $y$ be the vertex not on $C$. Then, by Lemma 6(i), the following two cases are possible: There is a vertex $x_i$ and an integer $a\geq 1$ such that $d(y,\{x_{i+1},x_{i+2},\ldots , x_{i+a}\})=0$, $x_{i-1}\rightarrow y \rightarrow x_i$ and the vertices $y$, $x_{i+a+1}$ are adjacent (Case I) or $d(y,\{x_{i+1},x_{i+2},\ldots , x_{i+a}\})=0$, $y\rightarrow \{x_{i}, x_{i+a+2}\}$, $x_{i+a+1}y\in D$ and the vertices $y$, $x_{i-1}$ are non-adjacent, where $a\in [1,n-6]$.

The proof will be by induction on $a$. We will first show that the theorem is true for $a=1$.

 \textbf{Case I.} $a=1$. Without loss of generality, we may assume that $x_{n-2}\rightarrow y\rightarrow x_{n-1}$, $x_2, y$ are adjacent and $y$ , $x_1$ are non-adjacent. Since the vertex $y$ is not adjacent with three consecutive vertices of $C$ (Lemma 9), it follows that $y, x_{n-3}$ also are non-adjacent. The condition (*) implies that $x_{n-2}x_1\notin D$, since $d(y)\leq n-2$ and $x_{n-2}y\in D$. 

We show that $x_1$ has a partner on $C[x_3,x_{n-2}]$. Assume that this is not the case. Then by Lemma 5(i) we have  $d(x_1)\geq n-1$, since $x_{n-2}x_1\notin D$ and $yx_1\notin D$. On the other hand, using Lemma 2(ii) and Remark (ii), we obtain  
$$
n-1\leq  d(x_1)=d(x_1,\{x_2,x_{n-1}\})+d(x_1,C[x_3,x_{n-2}])\leq n-2,
$$
which is a contradiction.

So, indeed $x_1$ has a partner on $C[x_3,x_{n-2}]$. Let the arc $x_kx_{k+1}\in C[x_3,x_{n-2}]$ be partner of $x_1$, i.e., $x_{k}\rightarrow x_1\rightarrow x_{k+1}$. Notice that $k\in [4,n-4]$ (by Lemma 11). If $yx_2\in D$, then $D(n,2)=[yx_{n-1};y x_2x_3 \ldots x_kx_1x_{k+1}\ldots x_{n-2}x_{n-1}]$, a contradiction. Assume therefore that $yx_2\notin D$. Then $x_2y$, $yx_3\in D$ and $y$, $x_4$ are non-adjacent, by Lemmas 9  and 6(i). This together with the condition (*) implies that $x_2x_4\notin D$, since $x_2y\in D$ and  $d(y)\leq n-2$. If $x_{n-2}x_2\in D$, then $C':=x_{n-2}x_2yx_3\ldots x_kx_1x_{k+1}\ldots $ $ x_{n-2}$ is a cycle of length $n-1$ missing the vertex $x_{n-1}$ for which $\{x_{n-2},y\}\rightarrow x_{n-1}$. Then $x_{n-1}x_2\in D$, by Lemmas 6(i) and 4, i.e., $x_{n-1}$ is adjacent with three consecutive vertices of $C'$, which is contrary to Lemma 9. Assume therefore that $x_{n-2}x_2\notin D$. Now we show that $x_2$ also has a partner on $C[x_3,x_{n-2}]$. Assume that this is not the case. Then, since $x_2x_{n-1}\notin D$ (by Lemma 11) and $x_2x_4\notin D$, using Lemma 2(iii) and Remark (ii) we obtain   
$$
 d(x_2)=d(x_2,\{y,x_1,x_3,x_{n-1}\})+d(x_2,C[x_4,x_{n-2}])\leq n-2.
$$
This together with $x_1\rightarrow \{x_2,x_{k+1}\}$ and the condition (*) implies that $x_2$, $x_{k+1}$ are adjacent. It is easy to see that $x_{k+1}x_2\in D$. By a similar argument, we conclude that $x_{n-2}x_2\in D$, which contradicts the fact that $x_{n-2}x_2\notin D$. Thus, $x_2$ also has a partner on $C[x_3,x_{n-2}]$. Therefore by Multi-Insertion Lemma there is a $(x_3,x_{n-1})$-path with vertex set $V(C)$, which together with the arcs $yx_{n-1}$ and $yx_3$ forms a Hamiltonian bypass. This completes the discussion of induction first step for ($a=1$) Case I.\\

Now we consider the induction first step for Case II.

 \textbf{Case II.} $a=1$. Without loss of generality, we may assume that $y\rightarrow \{x_3,x_{n-1}\}$, $x_2y\in D$ and $d(y,\{x_1,x_4,x_{n-2}\})=0$. By induction first step of Case I, we may assume that $y, x_5$ also are non-adjacent. This together with $d(y)\leq n-2$, $x_2y\in D$ and the condition (*) implies that 
$$d^+(x_2,\{x_4,x_5,x_{n-2}\})=0, \eqno (5)
$$
and hence, by Lemma 11, in particular, the vertices $x_2, x_4$ are non-adjacent. If $x_{n-2}x_1\in D$, then the cycle $C':=x_{n-2}x_1x_2yx_3\ldots x_{n-2}$ has length $n-1$ missing the vertex $x_{n-1}$ and $\{x_{n-2},y\}\rightarrow x_{n-1}\rightarrow x_1$, i.e., for the cycle $C'$ and vertex $x_{n-1}$ the considered induction first step of Case I holds. Assume therefore that $x_{n-2}x_1\notin D$. Then $x_1, x_{n-2}$ are non-adjacent (Lemma 11). It is not difficult to see that $x_1$ has a partner on $C[x_3,x_{n-3}]$. Indeed, for otherwise from Lemma 5(i) it follows that $d(x_{n-1})\geq n-1$ and hence by Lemma 2 and Remark (ii), we have 
$$
n-1\leq  d(x_1)=d(x_1,\{x_2,x_{n-1}\})+d(x_1,C[x_3,x_{n-3}])\leq n-2,
$$
which is a contradiction. Thus, indeed $x_1$ has a partner on $C[x_3,x_{n-3}]$. Let the arc $x_kx_{k+1}\in C[x_3,x_{n-3}]$ be  a partner of $x_1$. Note that $k\in [4,n-4]$ (by Lemma 11). Therefore neither  the vertex $x_2$ nor the arc $x_1x_2$   has a partner on $C[x_3,x_{n-1}]$ (for otherwise, by Multi-Insertion Lemma, there is an $(x_3,x_{n-1})$-path with vertex set $V(C)$, which together with the arcs $yx_3$ and $yx_{n-1}$ forms a Hamiltonian bypass). Recall that $a(x_2,x_4)=0$ and $x_2x_5\notin D$ (by (5). Now using Lemma 2(ii) and Remark (ii) we obtain  
$$
 d(x_2)=d(x_2,\{y,x_1,x_3\})+d(x_2,C[x_5,x_{n-1}])\leq n-2.
$$
 This together with $x_1 \rightarrow \{x_2,x_{k+1}\}$ and the condition (*) implies that $x_2$ and $x_{k+1}$ are adjacent. Then $x_{k+1}x_2\in D$ (if $x_2x_{k+1}\in D$, then the arc $x_1x_2$ has a partner on $C[x_3,x_{n-1}]$). By a similar argument, we conclude that $\{x_{n-2},x_{n-1}\}\rightarrow x_2$. Then $C':=x_{n-2}x_2yx_3x_4\ldots x_kx_1x_{k+1}\ldots x_{n-2}$ is a cycle of length $n-1$, which does not contain the vertex $x_{n-1}$ and $d(x_{n-1},\{x_{n-2},x_2,y\})=3$, a contradiction against Lemma 9 and hence, the discussion of case $a=1$ is completed.\\ 

  \textbf{The induction hypothesis}. Now we suppose that the  theorem is true if $D$ contains a cycle $C:=x_1x_2\ldots $ $ x_{n-1}x_1$ of length $n-1$ missing the vertex $y$ for which there is a vertex $x_i$ such that $d(y,\{x_{i+2},x_{i+3},\ldots ,$ $ x_{i+j}\})=0$ and 
(i) $x_i\rightarrow y\rightarrow x_{i+1}$ and the vertices $y, x_{i+j+1}$ are adjacent or
(ii) $y\rightarrow \{x_{i+1},x_{i+j+2}\}$ and $x_{i+j+1}y\in D$, where $2\leq j\leq a\leq n-6$.\\

Before dealing with Cases I and II, it is convenient to prove the following general Claim. 

 \textbf{Claim}.  Let $C:=x_1x_2\ldots x_{n-1}x_1$ be an arbitrary cycle of length $n-1$ in $D$ and let  $y$ be the vertex  not on $C$ and let $d(y,\{x_1,x_2,\ldots , x_a\})=0$, where $a\geq 2$.If 
(i) $x_{n-2}y,yx_{n-1}\in D$ and the vertices $y$ and $x_{a+1}$ are non-adjacent or 
(ii) $x_{a+1}y, yx_{a+2}$ and $yx_{n-1}\in D$, then  $x_{k-1}x_{k+1}\notin D$ for all $k\in [1,a]$.

\textbf{Proof of the claim}. Suppose, on the contrary, that $x_{k-1}x_{k+1}\in D$ for some $k\in [1,a]$, then $C':=x_{n-2}yx_{n-1}x_1\ldots x_{k-1}x_{k+1}\ldots$ $x_{n-2}$ or $C'':=x_{n-1}x_1\ldots x_{k-1}x_{k+1}\ldots x_{a+1}$ $yx_{a+2}\ldots x_{n-2}x_{n-1}$ is a cycle of length $n-1$ missing the vertex $x_k$ for (i) and (ii), respectively. Therefore, $d(x_k)\leq n-2$ (by Remark (iv)). By the induction hypothesis $x_k$ is not adjacent with vertices $x_{k+2},x_{k+3},\ldots ,x_{k+a},x_{k-2},x_{k-3},\ldots x_{k-a}$. In particular, $d^-(x_k,\{x_{k+1},x_{k+2},\ldots ,x_{a+1}\})=0$ and $d^-(x_k, C[x_{n-2},x_{k-2}])=0$ (it is easy to show that in both cases $x_{n-2}x_k\notin D$). Since $d(x_k)\leq n- 2$ and $x_{k-2}x_k\notin D$, by Lemma 5(i), the vertex $x_k$ has a partner on $C[x_{a+2},x_{n-2}]$, say the arc $x_jx_{j+1}\in C[x_{a+2},x_{n-2}]$ is a partner of $x_k$, i.e., $x_jx_k, x_kx_{j+1}\in D$. Therefore $x_{n-1}x_1\ldots x_{k-1}x_{k+1}\ldots x_ax_{a+1}\ldots x_jx_kx_{j+1}\ldots x_{n-2}x_{n-1}$ is a cycle of length $n-1$ missing the vertex $y$, for which $d(y, C[x_1,x_a]-\{x_k\})=0$ and $x_{n-2}y, yx_{n-1}\in D$, $x_{a+1}, y$ are non-adjacent or $yx_{n-1}, x_{a+1}y, yx_{a+2}\in D$ for (i) and (ii), respectively. Therefore, by the induction hypothesis $D$ contains a Hamiltonian bypass, a contradiction to our assumption. The claim is proved. \\

\textbf{Case I}. Without loss of generality, we may assume that $d(y,\{x_2,x_3,\ldots , x_{a+1}\})=0$, where $a\geq 2$, $x_{n-1}y, yx_1\in D$ and the vertices $y, x_{a+2}$ are non-adjacent. 

Notice, the condition (*) implies that for all $i\in [2,a+1]$,  $x_{n-1}x_i\notin D$, since $x_{n-1}y\in D$, $d(y)\leq n-2$ and the vertices $x_i, y$ are non-adjacent.

\textbf{Subcase I.1}. There are integers $k$ and $l$ with $1\leq l<k\leq a+2$ such that $x_kx_l\in D$. Without loss of generality, we assume that $k-l$ is as small as possible. From Remark (ii) and Lemma 11 it follows that $k-l\geq 3$. If every vertex $x_i\in C[x_{l+1},x_{k-1}]$ has a partner on the path $P:=x_kx_{k+1}\ldots x_{n-1}yx_1\ldots x_l$, then by Multi-Insertion Lemma there exists an $(x_k,x_l)$-Hamiltonian path, which together with the arc $x_kx_l$ forms a Hamiltonian bypass. Assume therefore that some vertex $x_i\in C[x_{l+1},x_{k-1}]$ has no partner on $P$. From the minimality of $k-l\geq 3$ and Claim 1 it follows that  $x_{i-2}\in C[x_l,x_k]$ and $a(x_i,x_{i-2})=0$ or $x_{i+2}\in C[x_l,x_k]$ and $a(x_i,x_{i+2})=0$. Therefore by the minimality of $k-l$ we have 
$$
d(x_i, C[x_l,x_k])\leq k-l-1. \eqno (6)
$$
Since $x_i$ has no partner on the path $C[x_{k+1},x_{n-1}]$, and if $l\geq 2$ also on $C[x_1,x_{l-1}]$, using Lemma 2 with the fact that $x_{n-1}x_i\notin D$ we obtain 
$$
d(x_i,C[x_{k+1},x_{n-1}]\leq n-k-1 \quad \hbox{and if} \quad l\geq 2, \quad \hbox{then} \quad d(x_i, C[x_1,x_{l-1}])\leq l.
$$
The last two inequalities together with (6) give: if $l\geq 2$, then $d(x_i)\leq n-2$, and if $l=1$, then $d(x_i)\leq n-3$. Thus, $d(x_i)\leq n-2$. In addition, Claim 1 and $x_{n-1}x_i\notin D$ imply that $x_{i-2}x_i\notin D$. Therefore, by Lemma 5(i), $x_i$ has a partner on $P$, which is contrary to our assumption.  

\textbf{Subcase I.2}. For any pair of  integers $k$ and $l$ with $1\leq l<k\leq a+2$,  $x_kx_l\notin D$. Then it is easy to see that   for each $x_i\in C[x_2,x_{a+1}])$,
$$
d(x_i, C[x_1,x_{a+2}])\leq a, \eqno (7) $$ 
 since $x_{i-2}\in C[x_1,x_{a+2}]$ and $a(x_i,x_{i-2})=0$ or $x_{i+2}\in C[x_1,x_{a+2}]$ and $a(x_i,x_{i+2})=0$. 

We first show that every vertex $x_{i}\in C[x_2,x_{a+1}]$ has a partner on $C[x_{a+3},x_{n-1}]$. Assume that this is not the case, i.e., some vertex $x_{i}\in C[x_2,x_{a+1}]$ has no partner on $C[x_{a+3},x_{n-1}]$. Then, since $x_{n-1}x_i\notin D$, by Lemma 2(ii) we have that $d(x_i,C[x_{a+3},x_{n-1}])\leq n-a-3$. This inequality together with (7) gives $d(x_i)\leq n-3$, a contradiction against Lemma 5(i), since $x_{i-2}x_i\notin D$. Thus each  vertex $x_{i}\in C[x_2,x_{a+1}]$ has a partner on $C[x_{a+3},x_{n-1}]$. Therefore, by Multi-Insertion Lemma there is an $(x_{a+3},x_{n-1})$-path, say $R$, with vertex set $V(C)-\{x_1,x_{a+2}\}$. If $yx_{a+2}\in D$, then $[yx_1; yx_{a+2}Rx_1]$ is a Hamiltonian bypass. Assume therefore that $yx_{a+2}\notin D$. Then $x_{a+2}y\in D$. By Lemma 6(i) and by the induction hypothesis, we have  $yx_{a+3}\in D$ and $d(y,\{x_{a+4},x_{a+5}\})=0$. This together with $x_{a+2}y\in D$, $d(y)\leq n-2$ and the condition (*) implies that 
$$
d^+(x_{a+2},,\{x_{a+4},x_{a+5}\})=0, \eqno (8) 
$$
in particular, by Lemma 11, $a(x_{a+2},a_{a+4})=0$. Since $yx_{a+3}\in D$ and each vertex $x_i\in C[x_2,x_{a+1}]$ has a partner on $C[x_{a+3},x_{n-1}]$, to show that $D$ contains a Hamiltonian bypass, by Multi-Insertion Lemma it suffices to prove that $x_{a+2}$ also has a partner on $C[x_{a+3},x_{n-1}]$. Assume that $x_{a+2}$  has no partner on $C[x_{a+3},x_{n-1}]$. Then, since the vertices $x_a$ and $x_{a+2}$ are non-adjacent (Claim 1 and Lemma 11), from Lemma 5(i) it follows that $d(x_{a+2})\geq n-1$. On the other hand, using (7), (8), $d(x_{a+2},\{x_a,x_{a+4}\})=0$ and Lemma 2, we obtain  
$$
n-1\leq d(x_{a+2})=d(x_{a+2},C[x_{1},x_{a+1}])+d(x_{a+2},\{y,x_{a+3})+d(x_{a+2},C[x_{a+5},x_{n-1}])\leq n-3,
$$
a contradiction. So, $x_{a+2}$ also has a partner on $C[x_{a+3},x_{n-1}]$ and the discussion of Case I is completed.

\textbf{Case II}. Without loss of generality, we assume that $d(y,\{x_1,x_2,\ldots , x_{a}\})=0$, where $a\geq 2$, $x_{a+1}y\in D$ and  $y\rightarrow \{x_{n-1} x_{a+2}\}$. 

By the considered Case I, without loss of generality, we may assume that $d(y,\{x_{n-2},x_{a+3}, x_{a+4}\})=0$. Since $x_{a+1}y\in D$, $d(y)\leq n-2$ and $d(y,\{x_1,x_2,\ldots , x_{a},x_{a+3},x_{a+4},x_{n-2}\})=0$, the condition (*) implies that  

$$d^+(x_{a+1},\{x_1,x_2,\ldots , x_{a},x_{a+3},x_{a+4},x_{n-2}\})=0. \eqno (9) $$

\textbf{Subcase II.1}. There are integers $k$ and $l$ with $1\leq l<k\leq a+1$ such that $x_kx_l\in D$. By (9), $k\not= a+1$. Without loss of generality, we assume that $k-l$ is as small as possible. By Remark (ii) and Lemma 11 we have $k-l\geq 3$.

We first show that each vertex of $C[x_{l+1},x_{k-1}]$ has a partner on the path $P:=x_kx_{k+1}\ldots x_{a+1}y$ $x_{a+2}\ldots x_{n-1}x_1\ldots x_l$. Assume that this is not the case and let $x_i\in C[x_{l+1},x_{k-1}]$ have no partner on $P$. Then, since $x_{i-2}x_i\notin D$ (Claim 1), from Lemma 5(i) and the minimality of $k-l$ it follows that $d(x_i)\geq n-1$. On the other hand, using the minimality of $k-l$ and the fact that  
$x_{i-2}\in C[x_l,x_k])$ and $a(x_i,x_{i-2})=0$ or  $x_{i+2}\in C[x_l,x_k]$ and $a(x_i,x_{i+2})=0$ we obtain 
 $$
d(x_i,C[x_{l},x_{k}])\leq k-l-1.$$
In addition,  by Lemma 2 and $x_{a+1}x_i\notin D$ we also have
 $$
d(x_i, C[x_{k+1},x_{a+1}])\leq a-k+1 \quad \hbox{and } \quad
d(x_i, C[x_{a+2},x_{l-1}])\leq n-a+l-2.$$  
Summing the last three inequalities gives $ d(x_i)\leq n-2$, which  contradicts that $d(x_i)\geq n-1$. Thus, indeed each  vertex $x_{i}\in C[x_{l+1},x_{k-1}]$ has a partner on $P$. Then, by Multi-Insertion Lemma, there is an $(x_k,x_l)$-Hamiltonian path , which together with the arc $x_kx_l$ forms a Hamiltonian bypass.

\textbf{Subcase II.2}. There are no $i$ and $j$ such that $1\leq i<j\leq a+1$ and  $x_jx_i\notin D$. If every vertex  $x_i\in C[x_1,x_{a+1}])$ has a partner on $C[x_{a+2},x_{n-1}]$, then by Multi-Insertion Lemma there is an  $(x_{a+2},x_{n-1})$-path, say $R$, with vertex set $V(C)$. Therefore $[yx_{n-1}; yR]$ is a Hamiltonian bypass. Assume therefore that there is a vertex $x_i\in C[x_1,x_{a+1}]$ which has no partner on $C[x_{a+2},x_{n-1}]$. 

Let $x_{i-2}x_i\notin D$, then from Lemma 5(i) it follows that $d(x_i)\geq n-1$. 

Assume first that $d(x_i,C[x_1,x_{a+1}])=a-1$. Using Lemma 2 we obtain that if $x_i\not=x_{a+1}$, then
 $$
n-1\leq d(x_{i})=d(x_{i},C[x_{1},x_{a+1}])+d(x_{i},C[x_{a+2},x_{n-1}])\leq n-2,$$ 
and, since  $x_{a+1}x_{a+3}\notin D$, if $x_i=x_{a+1}$, then
$$
n-1\leq d(x_{a+1})=d(x_{a+1},C[x_{1},x_{a+1}])+d(x_{a+1},\{y,x_{a+2})+d(x_{a+1},C[x_{a+3},x_{n-1}])\leq n-2,
$$
a contradiction.

Assume second that $d(x_i,C[x_1,x_{a+1}])=a$. Then from Claim 1 and Lemma 11 it follows that $a=2$, $x_i=x_2$, $d(x_2,\{x_1,x_3\})=2$ and $d(x_2,\{x_{n-1}\})=0$. Then

$$
n-1\leq d(x_{2})=d(x_{2},\{x_1,x_3\})+d(x_{2},C[x_{4},x_{n-2}])\leq n-2, $$
a contradiction.

Let now $x_{i-2}x_i\in D$. Then, by Claim 1, $x_i=x_1$ and $x_{n-2}x_1\in D$. We consider the cycle $C':=x_{n-3}x_{n-2}x_1x_2\ldots x_ax_{a+1}yx_{a+2}\ldots x_{n-3}$ of length $n-1$ missing the vertex $x_{n-1}$. Then $\{x_{n-2},y\}\rightarrow x_{n-1}$  and $x_{n-1}x_1\in D$, i.e., for the cycle $C'$ and the vertex $x_{n-1}$ Case I holds since $|\{x_2,x_3,\ldots ,x_{a+1} \}|=a$.  The discussion of Case II  is completed and with it the proof of the theorem is also completed. \fbox \\\\

\section { Concluding remarks }

The following two examples of digraphs show that if the minimal semi-degree of a digraph is equal to one, then the theorem is not true:

(i) Let $D(7)$ be a digraph with vertex set $\{x_1,x_2,\ldots ,x_6,y\}$ and let $x_1x_2\ldots x_6x_1$ be a cycle of length 6 in $D(7)$. Moreover, $N^+(y)=\{x_1,x_3,x_5\}$, $N^-(y)=\{x_2,x_4,x_6\}$, $x_1x_3, x_3x_5, x_5x_1\in D(7)$ and $D(7)$ has no other arcs. Note that $d^-(x_2)=d^-(x_4)=d^-(x_6)=1$ and $D(7)$ contains no dominated pair of non-adjacent vertices. It is not difficult to check that $D(7)$ contains no Hamiltonian bypass.

(ii) Let $D(n)$ be a digraph with vertex set $\{x_1,x_2,\ldots ,x_n\}$ and let $x_1x_2\ldots x_nx_1$ be a Hamiltonian cycle in $D(n)$. Moreover, $D(n)$ also contains the arcs $x_1x_3,x_3x_5, \ldots , x_{n-2}x_n$ (or  $x_1x_3,x_3x_5, \ldots , x_{n-3}x_{n-1},$ $ x_{n-1}x_1$ and $D(n)$ has no other arcs. Note that $D(n)$ contains no dominated pair of non-adjacent vertices, $d^-(x_2)=d^+(x_2)=1$. It is not difficult to check that $D(n)$ contains no Hamiltonian bypass. \\

We believe that Theorem 12 also is true if we require that the minimum in-degree at least two, instead of three .\\

In \cite{[2]} and \cite{[3]} Theorem 13 and Theorem 14 were proved, respectively.

\textbf{Theorem 13} (Bang-Jensen, Gutin, H. Li \cite{[2]}). Let $D$ be a strong digraph of order $n\geq 3$. Suppose that $min\{d^+(x)+d^-(y),d^-(x)+d^+(y)\}\geq n$ for any pair of non-adjacent vertices $x,y$ with a common out-neighbour or a common in-neighbour, then $D$ is Hamiltonian.

\textbf{Theorem 14} (Bang-Jensen, Guo, Yeo \cite{[3]}). Let $D$ be a strong digraph of order $n\geq 3$. Suppose that $d(x)+d(y)\geq 2n-1$ and  $min\{d^+(x)+d^-(y),d^-(x)+d^+(y)\}\geq n-1$ for any pair of non-adjacent vertices $x,y$ with a common out-neighbour or a common in-neighbour, then $D$ is Hamiltonian.\\
 
In \cite{[9]} and \cite{[10]} the following results were proved: 

\textbf{Theorem 15} ( \cite{[9]}). Let $D$ be a strong digraph of order $n\geq 4$ which is not a directed cycle. Suppose that $min\{d^+(x)+d^-(y),d^-(x)+d^+(y)\}\geq n$ for any pair of non-adjacent vertices $x,y$ with a common out-neighbour or a common in-neighbour. Then either $D$ contains a pre-Hamiltonian  cycle or $n$ is even and $D=K^*_{n/2,n/2}$.

\textbf{Theorem 16} ( \cite{[10]}). Let $D$ be a strong digraph of order $n\geq 4$ which is not a directed cycle. Suppose that $d(x)+d(y)\geq 2n-1$ and  $min\{d^+(x)+d^-(y),d^-(x)+d^+(y)\}\geq n-1$ for any pair of non-adjacent vertices $x,y$ with a common out-neighbour or a common in-neighbour. Then $D$ contains a pre-Hamiltonian  cycle or a cycle of length $n-2$. \\

In view of Theorems 13-16, we pose the following problem: 

\textbf{Problem}. Characterize those  digraphs which satisfy the condition of Theorem 13 or 14 but have no Hamiltonian bypass.

\textbf{Acknowledgment}. The autors will be grateful to the colleagues who will make any mathematical and grammatical comments.

\end{document}